\documentclass[12pt,twoside,reqno]{amsart}
\usepackage[colorlinks=true,citecolor=blue]{hyperref}
\usepackage{mathptmx, amsmath, amssymb, amsfonts, amsthm, mathptmx, enumerate, color}
\usepackage{graphicx}
\usepackage{epstopdf}

\newtheorem{theorem}{Theorem}[section]

\newtheorem{lemma}{Lemma}[section]
\newtheorem{proposition}{Proposition}[section]
\theoremstyle{definition}
\newtheorem{definition}{Definition}[section]

\newtheorem{remark}{Remark}[section]

\numberwithin{equation}{section}

\begin{document}
\setcounter{page}{1}

\vspace*{1.0cm}
\title[Contractibility  of the solution sets]
{Contractibility  of the solution sets   for set optimization problems}
\author[B. Chen, Y. Han]{Bin Chen$^{1}$, Yu Han$^{2,*}$}
\maketitle
\vspace*{-0.6cm}

\begin{center}
{\footnotesize {\it

$^1$School of Mathematical Sciences, Huaqiao University, Quanzhou, Fujian 362021,  China\\
$^2$ School of Statistics, Jiangxi University of Finance and Economics, Nanchang, Jiangxi 330013,  China

}}\end{center}

\vskip 4mm {\small\noindent {\bf Abstract.} This paper aims at investigating the contractibility  of the solution sets   for set optimization problems by   utilizing strictly  quasi cone-convexlikeness, which is an 
assumption weaker than	strictly cone-convexity, strictly cone-quasiconvexity and strictly naturally quasi cone-convexity. We establish the  contractibility  of  $l$-minimal,  $l$-weak minimal,  $u$-minimal  and  $u$-weak minimal     solution    sets   for set optimization problems by using the star-shape sets and the	nonlinear scalarizing functions for sets. Moreover, we also discuss the arcwise connectedness and 	the  contractibility of $p$-minimal  and  $p$-weak minimal     solution    sets   for set optimization problems by using the   scalarization technique. Finally, our main results are applied to the
contractibility  of the solution sets   for  vector optimization problems.

 }

\renewcommand{\thefootnote}{}
\footnotetext{ $^*$Corresponding author.
	
 E-mail addresses:chenbinmath@163.com (B. Chen), hanyumath@163.com (Y. Han).}

\section{Introduction}
It is known that investigating the properties of the solution set is one of the most important problems for optimization and related problems. Among many desirable properties
of the solution set, the connectedness is of considerable interest, since it provides
the possibility of continuously moving (transformations) from one solution to another.
Connectedness of the solution sets for vector optimization problems and vector equilibrium problems has been intensively studied in the literature, for instance, we refer the reader to \cite{Sun, LW, Gong1, Gong3, Gong4, Qiu,  HH4, HH2,  Han2}.
It is worth mentioning that the scalarization technique and  the scalarization results play important roles in investigating connectedness of solution sets of vector optimization problems and vector equilibrium problems (see, for example \cite{LW, Gong3,  Qiu, HH4, HH2,  Han2}).

Over the last two decades, set optimization problems have been extensively studied and applied in many fields such as optimal control problems, vector optimization problems, fuzzy optimization problems, image processing problems, viability theory, mathematical economics and differential inclusions (see \cite{KTZ15}).  Until now,    there are only three papers  \cite{HWH2019, Han2020, Han2022LT} considering connectedness of  solution sets  for set optimization problems. Han et al. \cite{HWH2019} established the arcwise connectedness of the sets of  minimal solutions and weak  minimal solutions for set optimization problems involving lower set less relation under the assumption of the strictly natural quasi cone-convexity.  Han \cite{Han2020}   obtained a scalarization result for weak $p$-minimal solution set for set optimization problems.   Moreover, he established connectedness of weak $p$-minimal solution set for set optimization problems by utilizing  the scalarization result.  Han \cite{Han2022LT} investigated connectedness of the sets of $l$-minimal approximate solutions and weak $l$-minimal approximate solutions for set optimization problems by using  the	nonlinear scalarizing functions for sets and the scalarization results.
However, to the best of our knowledge, the contractibility  of the solution sets   for set optimization problems has not been explored until now. Therefore, it is interesting and important to investigate the contractibility  of the solution sets   for set optimization problems.

It is known that the  nonlinear scalarizing functions play important roles in  set optimization problems.  Hern\'{a}ndez  and   Rodr\'{\i}guez-Mar\'{\i}n \cite{Hernandez1} introduced a nonlinear scalarizing function  with respect to lower set less relation. Moreover, they gave many interesting properties of the nonlinear scalarizing function.   Han and Huang  \cite{HHJOTA18}  established  the continuity and convexity of the nonlinear scalarizing function, which is introduced by Hern\'{a}ndez  and   Rodr\'{\i}guez-Mar\'{\i}n \cite{Hernandez1}.
As applications, they  investigated the upper semicontinuity and the lower semicontinuity of strongly approximate solution mappings to the parametric set optimization problems.
Very recently, Han  \cite{Han2019} established H\"{o}lder continuity of the nonlinear
scalarizing functions with respect to  lower set less relation and  upper set less relation, and discussed some properties of the nonlinear scalarizing function with respect to upper set less relation.
Moreover, he applied H\"{o}lder continuity of the nonlinear
scalarizing function to obtain Lipschitz continuity of strongly approximate solution mappings to the parametric set optimization problems.

At the perspective part of the paper \cite{HHJOTA18}, the authors   proposed that it may deserve study to explore the connectedness and the arcwise connectedness of the set of solutions for the set optimization problems by using   the nonlinear scalarizing function for sets.
Inspired by the perspective part of the paper \cite{HHJOTA18}, we make a new attempt to establish the  contractibility  of  $l$-minimal,  $l$-weak minimal,  $u$-minimal  and  $u$-weak minimal     solution    sets   for set optimization problems by using   the nonlinear scalarizing functions   introduced by Hern\'{a}ndez  and   Rodr\'{\i}guez-Mar\'{\i}n \cite{Hernandez1} and Han \cite{Han2019}.

On the other hand, we would like to point out that   the convexity and the generalized convexity of the objective mapping play important roles in investigating the stability of the    solution sets for  set optimization problems  (see, for example \cite{HH2017, HWH2020,  Han2020-11, KL2019, KL2020, ZH2020,  ZH2021}).   Very recently, Han and Zhang \cite{Han2022} introduced   the concepts of strictly quasi cone-convexlikeness and supremum mapping of the set-valued mappings, and established the upper semicontinuity and lower semicontinuity of the minimal solution mapping and weak minimal
solution mapping to parametric set optimization problems on Banach lattices by using strictly quasi cone-convexlikeness and the continuity of supremum mapping.  It is worth noting that strictly  quasi cone-convexlikeness is weaker than 
strictly cone-convexity, strictly cone-quasiconvexity and strictly naturally quasi cone-convexity.
In this paper,  we make a new attempt to investigate the contractibility  of the solution sets   for set optimization problems by   utilizing strictly  quasi cone-convexlikeness.

The rest of the paper is organized as follows. In Section 2, we present some necessary notations and lemmas.   In Section 3,  we derive the contractibility of  $l$-minimal and  $l$-weak minimal solution sets  for set optimization problems by using   the nonlinear scalarizing function   introduced by Hern\'{a}ndez  and   Rodr\'{\i}guez-Mar\'{\i}n \cite{Hernandez1}.   In Section 4, we establish the contractibility of  $u$-minimal and  $u$-weak minimal solution sets  for set optimization problems by using   the nonlinear scalarizing function   introduced by Han \cite{Han2019}. In Section 5, we investigate the arcwise connectedness and the contractibility of  $p$-minimal and  $p$-weak minimal solution sets  for set optimization problems by  utilizing   the scalarization technique. In Section 6,  we give an application of the main results  to obtain the
contractibility  of the solution sets  for vector optimization problems.

\section{Preliminaries}

Throughout this paper,   let $X$   be a normed vector space and let  $Y$ be a Banach space.  Assume that $C \subseteq Y$ is a convex, closed and pointed cone with nonempty interior, and $A$ and $B$ are two nonempty subsets of $Y$.    We denote by  $B_Y$  the closed unit ball in $Y$. Let ${{\wp _{0}}\left( Y \right)}$ be the family of all nonempty subsets of $Y$.
It is said that a nonempty set $A \subseteq Y$ is $C$-proper if $A + C \ne Y$, $C$-convex if $A + C$ is a convex set, $C$-closed if $A + C$ is a closed set, $C$-bounded if for each neighborhood $U$ of zero in $Y$ there is some positive number $t$ such that $A \subseteq tU + C$ and $C$-compact if any cover of $A$ of the form $\left\{ {{U_\alpha } + C:{U_\alpha }\;{\rm{are\;open}}} \right\}$ admits a finite subcover. Let ${{\wp _{0C}}\left( Y \right)}$ and ${\wp _{0 - C}}\left( Y \right)$ be the family of all nonempty $C$-proper and $-C$-proper subsets of $Y$, respectively.  Every $C$-compact set is $C$-bounded and $C$-closed (see \cite{Luc}). 

Let ${Y^*}$ be the topological dual space of $Y$ and ${C^*}$ be defined by
$${C^*} = \left\{ {f \in {Y^*}:f\left( c \right) \ge 0,\forall c \in C} \right\}.$$
A nonempty and convex subset $B$ of the convex cone $C$ is called a base of $C$, if $C = {\rm{cone}}\left( B \right)$ and $0 \notin {\rm{cl}}\left( B \right)$.
Let $e$ be a fixed point in ${\mathop{\rm int}} C$ and
$B_e^* = \left\{ {f \in {C^*}:f\left( e \right) = 1} \right\}$. It is clear that $B_e^*$ is a   base of ${C^*}$.

\begin{definition} \emph{(Kuroiwa, Tanaka, Ha (1997), \cite{KTH1997})}  Let $A,B \in {{\wp _{0}}\left( Y \right)}$. Then,
	$$A \le ^h B\mathop  {\Longleftrightarrow} \limits^{def} \left( {\bigcap\limits_{a \in A} {\left( {a + C} \right)} } \right) \cap B \ne \emptyset ,$$
	$$A{ \le ^l}B \mathop  {\Longleftrightarrow} \limits^{def} B \subseteq A + C,$$
	$$A{ \le ^u}B \mathop  {\Longleftrightarrow} \limits^{def} A \subseteq B - C.$$
	
\end{definition} 

Due to Proposition 2.3 of   \cite{KTY2008}, we know that the set relation $\le ^h$ is not  reflexive, but the set relation $\le ^h$ is  transitive.  Because of this, Han \cite{Han2020} introduced a set relation   $\le ^p$ as follows:

If $A = B$, then $A \le ^p B$;  if  $A \ne B$, then  $A \le ^p B \Leftrightarrow A \le ^h B$. It is easy to see that  $\le ^p$  is   reflexive and transitive.

The weak set relations "${ \ll ^l}$", "${ \ll ^u}$" and "${ \ll ^p}$", respectively, are defined by
$$A{ \ll ^l}B \mathop  {\Longleftrightarrow} \limits^{def} B \subseteq A + {\mathop{\rm int}} C,$$
$$A{ \ll ^u}B \mathop  {\Longleftrightarrow} \limits^{def} A \subseteq B - {\mathop{\rm int}} C,$$ and
$$A \ll ^p B\mathop  {\Longleftrightarrow} \limits^{def} \left( {\bigcap\limits_{a \in A} {\left( {a + C} \right)}  + {\mathop{\rm int}} C} \right) \cap B \ne \emptyset.$$

Recall that a nonempty set $A$ of a topological space is said to be arcwise connected, if  for every two points $x,y \in A$,  there exists a continuous mapping $\varphi :\left[ {0,1} \right] \to A$ such that $\varphi \left( 0 \right) = x$ and $\varphi \left( 1 \right) = y$.

\begin{definition}  A nonempty set $A$ of a topological space is said to be contractible, if  there exist ${x_0} \in A$ and a continuous mapping $H:A \times \left[ {0,1} \right] \to A$  such that $H\left( {x,1} \right) = x$ and 	
	$H\left( {x,0} \right) = {x_0}$ for any $x \in A$.	 
\end{definition}

\begin{definition} \cite{BS1986}  Let $A$ be a nonempty subset of $X$ and 
	$${\rm{star}}A: = \left\{ {x \in A:\lambda x + \left( {1 - \lambda } \right)y \in A,\;\forall \lambda  \in \left[ {0,1} \right],\;\forall y \in A} \right\}.$$
	We say that $A$	is starshaped when ${\rm{star}}A$ is nonempty.
\end{definition}

\begin{remark} \label{Wddfk}  It is clear that the assumption that $A$	is starshaped is   weaker than the assumption that $A$	is  convex.
  \end{remark}

Let $A \in {{\wp _{0}}\left( Y \right)}$ and $a \in A$. We say that $a$ is a minimal (resp. maximal) point of $A$ with respect to $C$ and we write $a \in {\rm{Min}}\left( A \right)$ (resp. $a \in {\rm{Max}}\left( A \right)$) if $\left( {A - a} \right) \cap \left( { - C} \right) = \left\{ 0 \right\}$ (resp. $\left( {A - a} \right) \cap C = \left\{ 0 \right\}$).
We say that $a$ is a weak minimal (resp. weak maximal) point of $A$ with respect to $C$ and we write $a \in {\rm{WMin}}\left( A \right)$ (resp. $a \in {\rm{WMax}}\left( A \right)$) if $\left( {A - a} \right) \cap \left( { - {\mathop{\rm int}} C} \right) = \emptyset $ (resp. $\left( {A - a} \right) \cap {\mathop{\rm int}} C = \emptyset $).

\begin{remark} \label{WMfk} It is clear that  ${\rm{Min}}\left( A \right) \subseteq {\rm{WMin}}\left( A \right)$ and ${\rm{Max}}\left( A \right) \subseteq {\rm{WMax}}\left( A \right)$.
	If $A$ is nonempty and $C$-compact, then ${\rm{Min}}\left( A \right) \ne \emptyset $ (see  \cite{Luc}). Similarly, 	if $A$ is nonempty and $-C$-compact, then ${\rm{Max}}\left( A \right) \ne \emptyset $.
\end{remark}

Let $F:X \to {2^Y}$  be a set-valued mapping and $K\subseteq X$ with $K\not=\emptyset$.
We   consider the following
set optimization problem:
\begin{center}(SOP) \quad  $\min F\left( x \right)$  \quad  subject to  \quad   $x \in K$.
\end{center}

\begin{definition} For each $s = l,u,p$,
	an element ${x_0} \in K$ is said to be
	\begin{itemize}
		\item[(i)]   $s$-minimal solution of  (SOP) if, for $x \in K$, $F\left( x \right) \le ^s F\left( {{x_0}} \right)$ implies $F\left( x_0 \right) {\le ^s} F\left( x \right) $.
		\item[(ii)]    weak $s$-minimal solution of  (SOP) if, for $x \in K$, $F\left( x \right) \ll ^s F\left( {{x_0}} \right)$  implies $F\left( x_0 \right) { \ll ^s}  F\left( x \right) $.
	\end{itemize}
\end{definition}
For each $s = l,u,p$, let ${E_s}\left( {F,K} \right)$ and  ${W_s}\left( {F,K} \right)$
denote the $s$-minimal solution set of (SOP) and  the weak $s$-minimal solution set of (SOP), respectively.

\begin{definition} \cite{HHJOTA18}  Let $D$ be a nonempty convex subset of $X$. A set-valued mapping $\Phi:X \to {2^Y}$ is said to be
	\begin{itemize}
		\item[(i)] naturally quasi $C$-convex on $D$ if, for any ${x_1},{x_2} \in D$   and for any $t \in \left[ {0,1} \right]$, there exists $\lambda  \in \left[ {0,1} \right]$ such that
		$$\lambda \Phi \left( {{x_1}} \right) + \left( {1 - \lambda} \right)\Phi \left( {{x_2}} \right) \subseteq \Phi \left( {t{x_1} + \left( {1 - t} \right){x_2}} \right) + C.$$
		\item[(ii)] strictly naturally quasi $C$-convex on $D$ if, for any ${x_1},{x_2} \in D$ with ${x_1} \ne {x_2}$  and for any $t \in \left( {0,1} \right)$, there exists $\lambda  \in \left[ {0,1} \right]$ such that
		$$\lambda \Phi \left( {{x_1}} \right) + \left( {1 - \lambda} \right)\Phi \left( {{x_2}} \right) \subseteq \Phi \left( {t{x_1} + \left( {1 - t} \right){x_2}} \right) + {\mathop{\rm int}} C.$$
	\end{itemize}
\end{definition}

\begin{definition} \label{ygst} \cite{Han2022}   Let $D$ be a nonempty    subset of $X$. For each $s = l,u,p$,  a set-valued mapping $\Phi:X \to {2^Y}$ is said to be
	\begin{itemize}
		\item[(i)]   $s$-$C$-convexlike  on $D$, if    for any ${x_1},{x_2} \in D$ and for any $t \in \left( {0,1} \right)$, there exists ${x_3} \in D$
		such that
		$$\Phi \left( {{x_3}} \right) \le ^s t\Phi \left( {{x_1}} \right) + \left( {1 - t} \right)\Phi \left( {{x_2}} \right).$$
		
		\item[(ii)]  strictly $s$-$C$-convexlike on $D$,  if for any ${x_1},{x_2} \in D$ with ${x_1} \ne {x_2}$ and for any $t \in \left( {0,1} \right)$, there exists ${x_3} \in D$ such that
		$$\Phi \left( {{x_3}} \right) \ll ^s t\Phi \left( {{x_1}} \right) + \left( {1 - t} \right)\Phi \left( {{x_2}} \right).$$
		
		\item[(iii)] strictly quasi $s$-$C$-convexlike on $D$, if for any ${x_1},{x_2} \in D$ with ${x_1} \ne {x_2}$,  there exist ${x_3} \in D$ and   $t_0 \in \left[ {0,1} \right]$ such that
		$$\Phi \left( {{x_3}} \right) \ll ^s {t_0} \Phi \left( {{x_1}} \right) + \left( {1 - {t_0}} \right)\Phi \left( {{x_2}} \right).$$
	\end{itemize}
\end{definition}

\begin{remark}   In Definiton \ref{ygst}, we do not need to require that  $D$ is convex or arcwise connected. \end{remark}

\begin{remark}  Strictly  quasi cone-convexlikeness is weaker than 
	strictly cone-convexity, strictly cone-quasiconvexity and strictly naturally quasi cone-convexity.  By the definition of strictly  quasi cone-convexlikeness, we can see that strictly  quasi cone-convexlikeness is a   very weak condition.
\end{remark}

\begin{definition}  \cite{GHTZ}    A set-valued mapping $\Phi:X \to {2^Y}$ is said to be
	\begin{itemize}
		\item[(i)]   upper semicontinuous (u.s.c.) at ${u_0} \in T$ if,
		for any neighborhood $V$ of $\Phi\left( {{u_0}} \right)$, there exists
		a neighborhood $U\left( {{u_0}} \right)$ of ${u_0}$ such that for every $u \in U\left( {{u_0}} \right)$, $\Phi\left( u \right) \subseteq  V$.
		\item[(ii)] $C$-upper semicontinuous ($C$-u.s.c.) at ${u_0} \in T$ if,
		for any neighborhood $V$ of $\Phi\left( {{u_0}} \right)$, there exists
		a neighborhood $U\left( {{u_0}} \right)$ of ${u_0}$ such that for every $u \in U\left( {{u_0}} \right)$, $\Phi \left( u \right) \subseteq V + C$.
		\item[(iii)] lower semicontinuous (l.s.c.) at ${u_0} \in T$ if, for any
		$x \in \Phi\left( {{u_0}} \right)$ and any neighborhood $V$ of
		$x$, there exists a neighborhood $U\left( {{u_0}} \right)$ of
		${u_0}$ such that for every $u \in U\left( {{u_0}} \right)$, $\Phi\left( u \right) \cap V \ne \emptyset$.
		\item[(iv)] $C$-lower semicontinuous ($C$-l.s.c.) at ${u_0} \in T$ if, for any
		$x \in \Phi\left( {{u_0}} \right)$ and any neighborhood $V$ of
		$x$, there exists a neighborhood $U\left( {{u_0}} \right)$ of
		${u_0}$ such that for every $u \in U\left( {{u_0}} \right)$, $\Phi \left( u \right) \cap \left( {V - C} \right) \ne \emptyset $.
	\end{itemize}
	We say that $\Phi$ is u.s.c., $C$-u.s.c., l.s.c.  and $C$-l.s.c. on $X$, if it is  u.s.c., $C$-u.s.c., l.s.c.  and $C$-l.s.c. at each point
	$x \in X$, respectively. We say that $\Phi$ is continuous on $X$ if it is both u.s.c. and l.s.c. on $X$. We say that $\Phi$ is $C$-continuous on $X$ if it is both $C$-u.s.c. and $C$-l.s.c. on $X$.
\end{definition}

\begin{remark} \label{C-czhuj}  It is clear that if  $\Phi$ is continuous at $x \in X$, then $\Phi$ is $C$-continuous at $x \in X$. The class of $C$-continuous	mappings is strictly larger than the class of continuous mappings.
\end{remark}

From Lemma 2.3 of \cite{HH2017} and Remark \ref{WMfk}, we can get the following two lemmas.
\begin{lemma} \label{Khrj}  If ${x_0} \in K$ and
	$F\left( {{x_0}} \right)$ is nonempty and $C$-compact, then
	${x_0} \in {W_l}\left( {F,K} \right)$ if and only if there  does not exist $y \in K$ satisfying $F\left( y \right){ \ll ^l}F\left( {{x_0}} \right)$.
\end{lemma}

\begin{lemma} \label{uuKhrj}  If ${x_0} \in K$ and
	$F\left( {{x_0}} \right)$ is nonempty and $-C$-compact, then
	${x_0} \in {W_u}\left( {F,K} \right)$ if and only if there  does not exist $y \in K$ satisfying $F\left( y \right){ \ll ^u}F\left( {{x_0}} \right)$.
\end{lemma}

\begin{lemma} \label{WExd}  Let $K$ be a nonempty   subset of  $X$. If  $F$ is strictly quasi  $l$-$C$-convexlike on $K$  with nonempty and $C$-compact values. Then	${W_l}\left( {F,K} \right) = {E_l}\left( {F,K} \right)$.
\end{lemma}
\emph{Proof} It suffices to show that ${W_l}\left( {F,K} \right) \subseteq {E_l}\left( {F,K} \right)$.
Let ${x_0} \in {W_l}\left( {F,K} \right)$ and $\bar{x} \in K$ with $F\left( {\bar x} \right){ \le ^l}F\left( {{x_0}} \right)$, 
which means that 
\begin{equation}\label{UgB1} F\left( {{x_0}} \right) \subseteq F\left( {{\bar x}} \right) + C.
\end{equation}
Suppose that $\bar{x} \ne {x_0}$. Since  $F$ is strictly quasi $l$-$C$-convexlike on $K$, there exist ${z_0} \in K$ and  $\lambda  \in \left[ {0,1} \right]$ such that
\begin{equation}\label{UgB2} \lambda F\left( {\bar x} \right) + \left( {1 - {\lambda}} \right)F\left( {{x_0}} \right) \subseteq F\left( {{z_0}} \right) + {\mathop{\rm int}} C.  \end{equation}
For any $u \in F\left( {{x_0}} \right)$,  it follows from  (\ref{UgB1})  that there exist $v \in F\left( {\bar x} \right)$
and  ${c_0} \in C$ such that $u = v + {c_0}$. Thanks to (\ref{UgB2}), we get $$\lambda v + \left( {1 - \lambda } \right)u \in F\left( {{z_0}} \right) + {\mathop{\rm int}} C.$$
Then
\begin{eqnarray*} u &=& \lambda u + \left( {1 - \lambda } \right)u = \lambda v + \left( {1 - \lambda } \right)u + \lambda {c_0} \\
	&\in& F\left( {{z_0}} \right) + {\mathop{\rm int}} C + C \subseteq F\left( {{z_0}} \right) + {\mathop{\rm int}} C.	
\end{eqnarray*}
By the arbitrariness of $u \in F\left( {{x_0}} \right)$, we have 
$$F\left( {{x_0}} \right) \subseteq F\left( {{z_0}} \right) + {\mathop{\rm int}} C.$$
This means that $F\left( {z_0} \right){ \ll ^l}F\left( {{x_0}} \right)$. We conclude from Lemma \ref{Khrj} that ${x_0} \notin {W_l}\left( {F,K} \right)$,
which contradicts  ${x_0} \in {W_l}\left( {F,K} \right)$. Hence, we have $\bar{x} = {x_0}$. It is easy to see that $F\left( {{x_0}} \right){ \le ^l}F\left( {\bar x} \right)$, and so ${x_0} \in {E_l}\left( {F,K} \right)$. This completes the proof.   \hfill$\Box$

\begin{remark}  \label{WExsqxd} In Lemma  \ref{WExd}, we do not need to assume that  $K$ is convex.  Since the class of $C$-compactness is strictly larger than the class of compactness and strictly  quasi cone convexlikeness is weaker than 	strictly naturally quasi cone-convexity,
 it is easy to see that Lemma  \ref{WExd}   improves Lemma 2.5 of \cite{HWH2019}.
\end{remark}

\begin{lemma} \label{WEuuxd}  Let $K$ be a nonempty   subset of  $X$. If  $F$ is strictly quasi  $u$-$C$-convexlike on $K$  with nonempty, $-C$-convex and $-C$-compact values. Then	${W_u}\left( {F,K} \right) = {E_u}\left( {F,K} \right)$.
\end{lemma}
\emph{Proof}  Let ${x_0} \in {W_l}\left( {F,K} \right)$ and $y_0 \in K$ with $F\left( {y_0} \right){ \le ^u}F\left( {{x_0}} \right)$, 
which  yields
\begin{equation}\label{UguB1} F\left( {{y_0}} \right) \subseteq F\left( {{ {x_0}}} \right) - C.
\end{equation}
Suppose that ${y_0} \ne {x_0}$. Noting that   $F$ is strictly quasi $u$-$C$-convexlike on $K$, there exist ${z_0} \in K$ and  $t  \in \left[ {0,1} \right]$ such that
\begin{equation}\label{UguB2} F\left( {{z_0}} \right) \subseteq tF\left( {{y_0}} \right) + \left( {1 - t} \right)F\left( {{x_0}} \right) - {\rm{int}}C. \end{equation}
For any $u \in F\left( {{z_0}} \right)$,  it follows from  (\ref{UguB2})  that there exist ${v_0} \in F\left( {{y_0}} \right)$, ${h_0} \in F\left( {{x_0}} \right)$
and  ${c_0} \in {\rm{int}}C$ such that
\begin{equation}\label{UguB3} u = t{v_0} + \left( {1 - t} \right){h_0} - {c_0}.
\end{equation}
Due to  (\ref{UguB1}), there exist ${s_0} \in F\left( {{x_0}} \right)$   and  ${\bar c} \in C$ such that 
\begin{equation}\label{UguB4}{v_0} = {s_0} - {\bar c}.
\end{equation}
Since $F\left( {{x_0}} \right)$ is $-C$-convex, we have 
\begin{equation}\label{UguB5} t{s_0} + \left( {1 - t} \right){h_0} \in F\left( {{x_0}} \right) - C.
\end{equation}
Thanks to (\ref{UguB3}),   (\ref{UguB4}) and (\ref{UguB5}),    we get
$$u = t\left( {{s_0} - {\bar c}} \right) + \left( {1 - t} \right){h_0} - {c_0} \in F\left( {{x_0}} \right) - C - t{\bar c} - {\mathop{\rm int}} C \subseteq F\left( {{x_0}} \right) - {\mathop{\rm int}} C.$$
Due to the arbitrariness of $u \in F\left( {{z_0}} \right)$, we have 
$$F\left( {{z_0}} \right) \subseteq F\left( {{x_0}} \right) - {\mathop{\rm int}} C,$$
 which implies $F\left( {z_0} \right){ \ll ^u}F\left( {{x_0}} \right)$. Combining this with  Lemma \ref{uuKhrj}, we get ${x_0} \notin {W_u}\left( {F,K} \right)$.
which contradicts  ${x_0} \in {W_u}\left( {F,K} \right)$. Thus, we have ${y_0} = {x_0}$. This  implies  that $F\left( {{x_0}} \right){ \le ^u}F\left( {y_0} \right)$, and so ${x_0} \in {E_u}\left( {F,K} \right)$. Therefore, we have ${W_u}\left( {F,K} \right) \subseteq {E_u}\left( {F,K} \right)$. This together with   Remark 2.1 of \cite{HH2017} implies that ${W_u}\left( {F,K} \right) = {E_u}\left( {F,K} \right)$.
This completes the proof.   \hfill$\Box$

\begin{remark} In Lemma  \ref{WEuuxd}, we do not need to assume that  $K$ is convex. Moreover, noting that strictly  quasi cone-convexlikeness is weaker than 	strictly  cone-convexity, we can  see that Lemma  \ref{WEuuxd}   improves Proposition 2.2 of \cite{HH2017}.
\end{remark}

The following lemma is well known.
\begin{lemma} \label{HltJB} Let $\Upsilon $ and $\Omega $ be two topological spaces. Assume that $A$ is a nonempty and arcwise connected subset of $\Upsilon $ and $f:A \to \Omega$ is a continuous  single-valued   mapping. Then $f\left( A \right): = \bigcup\limits_{x \in A} {\left\{ {f\left( x \right)} \right\}} $	 is arcwise connected.
\end{lemma}

\section{Contractibility of  l-minimal and  l-weak minimal solution sets}
In this section,   we establish the contractibility of  $l$-minimal and  $l$-weak minimal solution sets  for set optimization problems.

\begin{definition}  \cite{GHTZ}  Let $a \in Y$. A function ${\phi _{e,a}}:Y \to \mathbb{R}$ is said to be the Gerstewizt's function if
	$${\phi _{e,a}}\left( y \right) = \min \left\{ {t \in \mathbb{R}:y \in te + a + C} \right\},\;\;\forall y \in Y.$$
\end{definition}

Replacing $a$ by a nonempty set $A \subseteq Y$, we obtain the function ${\phi _{e,A}}:Y \to \mathbb{R} \cup \left\{ { - \infty } \right\}$ as follows
$${\phi _{e,A}}\left( y \right) = \inf \left\{ {t \in \mathbb{R}:y \in te + A + C} \right\},\;\; \forall y \in Y.$$

\begin{definition}  \cite{Hernandez1}  Define the function ${G_e}\left( { \cdot , \cdot } \right):{{\wp _{0C}}\left( Y \right)} \times {{\wp _{0C}}\left( Y \right)} \to \mathbb{R} \cup \left\{ { - \infty , + \infty } \right\}$ by
	$${G_e}\left( {A,B} \right) = \mathop {\sup }\limits_{b \in B} \left\{ {{\phi _{e,A}}\left( b \right)} \right\},\;\; \forall \left( {A,B} \right) \in {{\wp _{0C}}\left( Y \right)} \times {{\wp _{0C}}\left( Y \right)}.$$
\end{definition}

We  define the function $\xi :K \times K \to \mathbb{R} \cup \left\{ { - \infty , + \infty } \right\}$ by
$$\xi \left( {x,y} \right) = {G_e}\left( {F\left( x \right),F\left( y \right)} \right), \quad \forall \left( {x,y} \right) \in K \times K.$$

It is easy to get the following lemma.

\begin{lemma} \label{phixz} Let $a,{a_1},{a_2} \in Y$, $y,{y_1},{y_2} \in Y$ and $\lambda  \ge 0$. Then:
	\begin{itemize}
		\item[(i)] ${\phi _{e,\lambda a}}\left( {\lambda y} \right) = \lambda {\phi _{e,a}}\left( y \right)$.
		\item[(ii)] ${\phi _{e,{a_1} + {a_2}}}\left( {{y_1} + {y_2}} \right) \le {\phi _{e,{a_1}}}\left( {{y_1}} \right) + {\phi _{e,{a_2}}}\left( {{y_2}} \right).$
		\item[(iii)] if ${a_2} - {a_1} \in C$, then ${\phi _{e,{a_1}}}\left( y \right) \le {\phi _{e,{a_2}}}\left( y \right)$.
		\item[(iv)]  if ${a_2} - {a_1} \in {\mathop{\rm int}} C$, then ${\phi _{e,{a_1}}}\left( y \right) < {\phi _{e,{a_2}}}\left( y \right)$.
	\end{itemize}
\end{lemma}

\begin{lemma} \label{tuxing2} Let $K$ be  a nonempty subset of $X$ and ${y _0} \in K$. Assume that   $F$ is 	strictly quasi $l$-$C$-convexlike on $K$  with nonempty  $C$-compact values. Then $\xi \left( {\cdot, {y _0} } \right)$ is a 	strictly quasi  convexlike function on $K$, i.e., for any ${x_1},{x_2} \in K$ with ${x_1} \ne {x_2}$,  there exist ${x_3} \in K$ and   $\lambda \in \left[ {0,1} \right]$ such that
	$$\xi \left( {{x_3},{y_0}} \right) < \lambda \xi \left( {{x_1},{y_0}} \right) + \left( {1 - \lambda } \right)\xi \left( {{x_2},{y_0}} \right).$$
\end{lemma}
{\it Proof}   Let ${x_1},{x_2} \in K$ with ${x_1} \ne {x_2}$.   Since  $F$ is 	strictly quasi $l$-$C$-convexlike on $K$, there exist ${x_3} \in K$ and   $t_0 \in \left[ {0,1} \right]$ such that
\begin{equation}\label{kkB1} {t_0}F\left( {{x_1}} \right) + \left( {1 - {t_0}} \right)F\left( {{x_2}} \right) \subseteq F\left( {{x_3}} \right) + {\mathop{\rm int}} C. 
\end{equation}
It follows from $C$-compactness of  $F\left( {{x_3}} \right)$  that $F\left( {{x_3}} \right)$
is $C$-closed. Due to $C$-compactness of $F\left( {{y _0}} \right)$ and  Proposition 3.4 of \cite{Hernandez1}, there exists  ${z_0} \in F\left( {{y_0}} \right)$ such that
\begin{equation}\label{kkB2}
{G_e}\left( {F\left( {{x_3}} \right),F\left( {{y_0}} \right)} \right) = \mathop {\max }\limits_{z \in F\left( {{y_0}} \right)} \left\{ {{\phi _{e,F\left( {{x_3}} \right)}}\left( z \right)} \right\}  
= {\phi _{e,F\left( {{x_3}} \right)}}\left( {{z_0}} \right).
\end{equation}
Since $F\left( {{x_1}} \right)$ and $F\left( {{x_2}} \right)$ are  $C$-compact, it follows from Proposition 2.18 (i) of \cite{Hernandez1} that there are  ${u_1} \in F\left( {{x_1}} \right)$ and   ${u_2} \in F\left( {{x_2}} \right)$  such that
\begin{equation}\label{kkB3} {\phi _{e,F\left( {{x_i}} \right)}}\left( {{z_0}} \right) = {\phi _{e,{u_i}}}\left( {{z_0}} \right),  \quad i = 1,2. 
\end{equation}
Due to ${u_1} \in F\left( {{x_1}} \right)$,   ${u_2} \in F\left( {{x_2}} \right)$ and (\ref{kkB1}), there exists  ${{u_t} \in F\left( {{x_3}} \right)}$ such that ${t_0} {u_1} + \left( {1 - t_0 } \right){u_2} - {u_t} \in {\mathop{\rm int}} C$.
Combining this with  Lemma  \ref{phixz} and (\ref{kkB3}), we get
\begin{eqnarray}\label{kkB4}
{\phi _{e,{u_t}}}\left( {{z_0}} \right) &<& {\phi _{e,{t_0} {u_1} + \left( {1 - {t_0} } \right){u_{_2}}}}\left( {{z_0}} \right) 
= {\phi _{e,{t_0} {u_1} + \left( {1 - {t_0} } \right){u_{_2}}}}\left( {{t_0} {z_0} + \left( {1 - {t_0} } \right){z_0}} \right) \nonumber \\
& \le& {t_0} {\phi _{e,{u_1}}}\left( {{z_0}} \right) + \left( {1 - {t_0} } \right){\phi _{e,{u_2}}}\left( {{z_0}} \right)  
= {t_0} {\phi _{e,F\left( {{x_1}} \right)}}\left( {{z_0}} \right) +  \left( {1 - {t_0} } \right){\phi _{e,F\left( {{x_2}} \right)}}\left( {{z_0}} \right).
\end{eqnarray}

We conclude from  (\ref{kkB2})  and  (\ref{kkB4})  that
\begin{eqnarray*}   {G_e}\left( {F\left( {{x_3}} \right),F\left( {{y_0}} \right)} \right) 
	&=& {\phi _{e,F\left( {{x_3}} \right)}}\left( {{z_0}} \right)  \le {\phi _{e,{u_t}}}\left( {{z_0}} \right) \\
	& <& {t_0} {\phi _{e,F\left( {{x_1}} \right)}}\left( {{z_0}} \right) + \left( {1 - {t_0} } \right){\phi _{e,F\left( {{x_2}} \right)}}\left( {{z_0}} \right)\\
	&\le& {t_0} {G_e}\left( {F\left( {{x_1}} \right),F\left( {{y_0}} \right)} \right) + 
	\left( {1 - {t_0} } \right){G_e}\left( {F\left( {{x_2}} \right),F\left( {{y_0}} \right)} \right).
\end{eqnarray*}
This means that
$$\xi \left( {{x_3},{y_0}} \right) < {t_0}\xi \left( {{x_1},{y_0}} \right) + \left( {1 - {t_0}} \right)\xi \left( {{x_2},{y_0}} \right).$$
This completes the proof.  \hfill$\Box$

\begin{remark} \label{dhkg}  	Similar to	Remark \ref{WExsqxd}, we can see that Lemma  \ref{tuxing2}   improves Theorem 4.4  (ii) of \cite{HHJOTA18}.
\end{remark}

From Theorem 5.1 of \cite{Han2019}, it is easy to get the following lemma.
\begin{lemma} \label{lianxu} If $F$ is $C$-continuous  on $K$ with nonempty  and $C$-compact  values, then   $\xi \left( { \cdot , \cdot } \right)$ is continuous on $K \times K$.
\end{lemma}

\begin{theorem} \label{hldtx}  Let $K$ be a nonempty  and compact subset of $X$. Assume that  $K$ is    starshaped   and $F$ is  $C$-continuous and strictly quasi $l$-$C$-convexlike  on $K$  with nonempty  $C$-compact values. Then ${W_l}\left( {F,K} \right)$ is contractible.  Moreover, ${E_l}\left( {F,K} \right)$ is  contractible.
\end{theorem}
{\it Proof} Since  $K$ is    starshaped, we have ${\rm{star}}K \ne \emptyset $. Let ${x_0} \in {\rm{star}}K$.
For any $x \in {W_l}\left( {F,K} \right)$ and for any $\lambda  \in \left[ {0,1} \right]$, let $\eta \left( {x,\lambda } \right) = \lambda x + \left( {1 - \lambda } \right){x_0}$. Due to ${x_0} \in {\rm{star}}K$, we have $\eta \left( {x,\lambda } \right) \in K$. We   define a set-valued mapping $H:{W_l}\left( {F,K} \right) \times \left[ {0,1} \right] \to {2^K}$  as follows
$$H(x,\lambda ) = \left\{ {v \in K:\xi \left( {y,\eta \left( {x,\lambda } \right)} \right) \ge \xi \left( {v,\eta \left( {x,\lambda } \right)} \right),\>\forall y \in K} \right\}.$$
It follows from Lemma  \ref{lianxu}  that $\xi \left( { \cdot , \cdot } \right)$ is continuous on $K \times K$. This together with  compactness of $K$  implies that $H\left( x,\lambda \right)$ is nonempty  for any $x \in {W_l}\left( {F,K} \right)$ and for any $\lambda  \in \left[ {0,1} \right]$.

We prove that $H\left( {x,\lambda } \right) \subseteq {W_l}\left( {F,K} \right)$ for any $x \in {W_l}\left( {F,K} \right)$ and for any $\lambda  \in \left[ {0,1} \right]$. In fact, let ${v_0} \in H\left( {{x,\lambda}} \right)$. Then
\begin{equation}\label{LThX1} {G_e}\left( {F(y),F\left( {{\eta \left( {x,\lambda } \right)}} \right)} \right) \ge {G_e}\left( {F\left( {{v_0}} \right),F\left( {{\eta \left( {x,\lambda } \right)}} \right)} \right),\quad \forall y \in K.
\end{equation}
Suppose that ${v_0} \notin {W_l}(F,K)$. Thanks to Lemma \ref{Khrj},   there is ${y_0} \in K$ such that $F\left( {{y_0}} \right){ \ll ^l}F\left( {{v_0}} \right)$.
It follows from Theorem 3.9 (ii) of \cite{Hernandez1} that
$${G_e}\left( {F\left( {{y_0}} \right),F\left( {\eta \left( {x,\lambda } \right)} \right)} \right) < {G_e}\left( {F\left( {{v_0}} \right),F\left( {\eta \left( {x,\lambda } \right)} \right)} \right),$$
which contradicts (\ref{LThX1}). Hence ${v_0} \in {W_l}(F,K)$ and so  $H\left( {x,\lambda } \right) \subseteq {W_l}\left( {F,K} \right)$.

We show that $H\left( x,\lambda \right)$ is a singleton   for any $x \in {W_l}\left( {F,K} \right)$ and for any $\lambda  \in \left[ {0,1} \right]$.
In fact, suppose that $H\left( x,\lambda \right)$ is not a singleton. Then  there exist  ${z_1},{z_2} \in H\left( x,\lambda \right)$  such that ${z_1} \ne {z_2}$.
It follows from Lemma  \ref{tuxing2} that $\xi \left( {\cdot, {\eta \left( {x,\lambda } \right)} } \right)$ is a 	strictly quasi  convexlike function on $K$. This means that
there exist ${z_3} \in K$ and   $\lambda \in \left[ {0,1} \right]$ such that
\begin{equation}\label{LThX2} \xi \left( {{z_3},\eta \left( {x,\lambda } \right)} \right) < \lambda \xi \left( {{z_1},\eta \left( {x,\lambda } \right)} \right) + \left( {1 - \lambda } \right)\xi \left( {{z_2},\eta \left( {x,\lambda } \right)} \right).   \end{equation}
Due to ${z_1},{z_2} \in H\left( x,\lambda \right)$ and ${z_3} \in K$, one has 
$$\xi \left( {{z_3},\eta \left( {x,\lambda } \right)} \right) \ge \xi \left( {{z_i},\eta \left( {x,\lambda } \right)} \right),\quad i = 1,2,$$
and so
$$\xi \left( {{z_3},\eta \left( {x,\lambda } \right)} \right) \ge \lambda \xi \left( {{z_1},\eta \left( {x,\lambda } \right)} \right) + \left( {1 - \lambda } \right)\xi \left( {{z_2},\eta \left( {x,\lambda } \right)} \right),$$
which contradicts (\ref{LThX2}). Therefore, $H\left( x,\lambda \right)$ is a singleton.

We can see that $H$ is single-valued  mapping from ${W_l}\left( {F,K} \right) \times \left[ {0,1} \right]$ to ${W_l}\left( {F,K} \right)$. Next, we prove that $H$ is continuous on ${W_l}\left( {F,K} \right) \times \left[ {0,1} \right]$. Suppose on the contrary that  $H$ is not continuous at $\left( {{u_0},{\lambda _0}} \right) \in {W_l}\left( {F,K} \right) \times \left[ {0,1} \right]$. Then there exist a neighborhood $W_0$ of $H {\left( {{u_0},{\lambda _0}} \right)}$
and a sequence $\left\{ {\left( {{u_n},{\lambda _n}} \right)} \right\} \in {W_l}\left( {F,K} \right) \times \left[ {0,1} \right]$
with $\left( {{u_n},{\lambda _n}} \right) \to \left( {{u_0},{\lambda _0}} \right)$ such that 
\begin{equation}\label{LThX3}
{v_n}: = H\left( {{u_n},{\lambda _n}} \right) \notin {W_0}. 
\end{equation}
Since $K$ is   compact,  without loss of generality,  we assume that ${v_n} \to {v_0} \in K$. For any   $y \in K$,  it follows from  ${v_n} = H\left( {{u_n},{\lambda _n}} \right)$  that
\begin{equation}\label{LThX4}  \xi \left( {y,\eta \left( {{u_n},{\lambda _n}} \right)} \right) \ge \xi \left( {{v_n},\eta \left( {{u_n},{\lambda _n}} \right)} \right).
\end{equation}
In view of Lemma  \ref{lianxu}, we get that  $\xi \left( { \cdot , \cdot } \right)$ is continuous on $K \times K$. This together with  (\ref{LThX4}) and $\eta \left( {{u_n},{\lambda _n}} \right) \to \eta \left( {{u_0},{\lambda _0}} \right)$ implies that
$$\xi \left( {y,\eta \left( {{u_0},{\lambda _0}} \right)} \right) \ge \xi \left( {{v_0},\eta \left( {{u_0},{\lambda _0}} \right)} \right),\quad \forall y \in K,$$
which yields ${v_0} = H\left( {{u_0},{\lambda _0}} \right)$. Hence, ${v_n} \to {v_0} = {H\left( {{u_0},{\lambda _0}} \right)} \subseteq {W_0}$ and so ${v_n} \in {W_0}$ for  $n$  large enough, which contradicts (\ref{LThX3}).

It is clear that $H(x,0) = H({x_0},0)$ for any   $x \in {W_l}\left( {F,K} \right)$.  Finally, it suffices to show that $H(x,1) = x$ for any   $x \in {W_l}\left( {F,K} \right)$.  For any $x \in {W_l}(F,K)$. due to Lemma \ref{Khrj}, we obtain that  for any $y \in K$, $F\left( y \right){ \ll ^l}F\left( {{x}} \right)$ does not hold. Thanks to Corollary 3.11 (ii) of  \cite{Hernandez1}, we get
\begin{equation}\label{LThX5} {G_e}\left( {F\left( y \right),F\left( {x} \right)} \right) \ge 0,\quad \forall y \in K.
\end{equation}
It follows from Theorem 3.10 (i) of \cite{Hernandez1} that ${G_e}\left( {F\left( {x} \right),F\left( {x} \right)} \right) = 0$.  Combining this with   (\ref{LThX5}), we have
\begin{equation}\label{LThX6}{G_e}\left( {F(y),F\left( {x} \right)} \right) \ge {G_e}\left( {F\left( {x} \right),F\left( {x} \right)} \right),\quad \forall y \in K.\end{equation}
We conclude from  (\ref{LThX6})  and  $\eta \left( {x,1} \right) = x$  that
$$\xi \left( {y,\eta \left( {x,1} \right)} \right) = \xi \left( {y,x} \right) \ge \xi \left( {x,x} \right) = \xi \left( {x,\eta \left( {x,1} \right)} \right),  \quad \forall y \in K,$$
which yields that  $x = H(x,1)$  for any   $x \in {W_l}\left( {F,K} \right)$. Therefore, ${W_l}\left( {F,K} \right)$ is contractible. This together with   Lemma \ref{WExd} implies that  ${E_l}\left( {F,K} \right)$ is contractible. This completes the proof.   \hfill$\Box$

\section{Contractibility of  u-minimal and  u-weak minimal solution sets}

In this section, we investigate the contractibility of  $u$-minimal and  $u$-weak minimal solution sets for set optimization problems.
\begin{definition} Let $w \in {\mathop{\rm int}} C$ and $b \in Y$. The  nonlinear scalarizing function ${\varphi _{w,b}}:Y \to \mathbb{R}$ is defined by
	$${\varphi _{w,b}}\left( y \right) = \min \left\{ {t \in \mathbb{R}:y \in tw + b - C} \right\},\;\;\forall y \in Y.$$
\end{definition}

Replacing $b$ by a nonempty set $B \subseteq Y$, we obtain the function ${\varphi _{w,B}}:Y \to \mathbb{R} \cup \left\{ { - \infty } \right\}$ as follows
$${\varphi _{w,B}}\left( y \right) = \inf \left\{ {t \in \mathbb{R}:y \in tw + B - C} \right\},\;\; \forall y \in Y.$$
Han \cite{Han2019} introduced  the  nonlinear scalarizing function ${H_w}\left( { \cdot , \cdot } \right):{\wp _{0 - C}}\left( Y \right) \times {\wp _{0 - C}}\left( Y \right) \to \mathbb{R} \cup \left\{ { - \infty , + \infty } \right\}$ for $u$-type less order relation defined by
$${H_w}\left( {A,B} \right) = \mathop {\sup }\limits_{a \in A} \left\{ {{\varphi _{w,B}}\left( a \right)} \right\},\;\;{\rm{for}}\left( {A,B} \right) \in {\wp _{0 - C}}\left( Y \right) \times {\wp _{0 - C}}\left( Y \right).$$

We  define the function $\gamma  :K \times K \to \mathbb{R} \cup \left\{ { - \infty , + \infty } \right\}$ by
$$\gamma \left( {x,y} \right) = {H_w}\left( {F\left( x \right),F\left( y \right)} \right), \quad \forall \left( {x,y} \right) \in K \times K.$$

\begin{proposition}  \label{Bxsdu} \cite{Han2019} Let $r \in \mathbb{R}$ and $B \in {{\wp _{0}}\left( Y \right)}$. Then for any $y \in Y$, we have:
	\begin{itemize}
		\item[(i)] ${\varphi _{w,B}}\left( y \right) < r \Leftrightarrow y \in rw + B - {\mathop{\rm int}} C$.
		\item[(ii)] ${\varphi _{w,B}}\left( y \right) \le r \Leftrightarrow y \in rw + {\rm{cl}}\left( {B - C} \right)$.
		\item[(iii)] ${\varphi _{w,B}}\left( y \right) \ge r \Leftrightarrow y \notin rw + B - {\rm{int}}C$.
		
		\item[(iv)] ${\varphi _{w,B}}\left( y \right) > r \Leftrightarrow y \notin rw + {\rm{cl}}\left( {B - C} \right).$
		\item[(v)] ${\varphi _{w,B}}\left( y \right) = r \Leftrightarrow y \in rw + \partial \left( {B - C} \right).$
	\end{itemize}
\end{proposition}

From Proposition  \ref{Bxsdu} (i), it is easy to get the following proposition.
\begin{proposition}  \label{BxYGsdu}  Let $r \in \mathbb{R}$ and $\left( {A,B} \right) \in {\wp _{0 - C}}\left( Y \right) \times {\wp _{0 - C}}\left( Y \right)$.  If ${H_w}\left( {A,B} \right) < r$, then $A \subseteq rw + B - {\mathop{\rm int}} C$. In particular,  we have ${H_w}\left( {A,B} \right) < 0 \Rightarrow A{ \le ^u}B$.	
\end{proposition}

From Proposition  \ref{Bxsdu} (ii), it is easy to get the following proposition.

\begin{proposition}  \label{Bxugsdu}  Let $r \in \mathbb{R}$ and $\left( {A,B} \right) \in {\wp _{0 - C}}\left( Y \right) \times {\wp _{0 - C}}\left( Y \right)$. Then
	$${H_w}\left( {A,B} \right) \le r \Leftrightarrow A \subseteq rw + {\rm{cl}}\left( {B - C} \right).$$
\end{proposition}

\begin{proposition}  \label{Bxxxsdu}   Let $y \in Y$ and $A,B \in {\wp _{0 - C}}\left( Y \right)$.
	 \begin{itemize}
		\item[(i)] If $A{ \le ^u}B$, then ${\varphi _{w,B}}\left( y \right) \le {\varphi _{w,A}}\left( y \right)$.			
		\item[(ii)] If $A{ \ll ^u}B$ and $A$ is $-C$-closed, then ${\varphi _{w,B}}\left( y \right) < {\varphi _{w,A}}\left( y \right)$.			
	\end{itemize}
\end{proposition}
\emph{Proof} (i). Let $r = {\varphi _{w,A}}\left( y \right)$. For any $\varepsilon  > 0$, we have ${\varphi _{w,A}}\left( y \right) < r + \varepsilon $. It follows from Proposition  \ref{Bxsdu} (i) that   
\begin{equation}\label{DxsBa1} y \in \left( {r + \varepsilon } \right)w + A - {\mathop{\rm int}} C.
\end{equation}
Due to $A{ \le ^u}B$, we have $A \subseteq B - C$. This together with  (\ref{DxsBa1})   implies that 
$$y \in \left( {r + \varepsilon } \right)w + B - C - {\mathop{\rm int}} C \subseteq \left( {r + \varepsilon } \right)w + B - C,$$
which yields ${\varphi _{w,B}}\left( y \right) \le r + \varepsilon $. By the arbitrariness of  $\varepsilon  > 0$, one has ${\varphi _{w,B}}\left( y \right) \le r = {\varphi _{w,A}}\left( y \right)$.

(ii) Let $\alpha  = {\varphi _{w,A}}\left( y \right)$. By  Proposition  \ref{Bxsdu} (ii) and the $-C$-closedness of $A$, we have 
\begin{equation}\label{DxsBa2} y \in \alpha w + {\rm{cl}}\left( {A - C} \right) \subseteq \alpha w + A - C.
\end{equation}
 It follows from  $A{ \ll ^u}B$  that $A \subseteq B - {\mathop{\rm int}} C$. Combining this with  (\ref{DxsBa2}), we get
 $$y \in \alpha w + B - {\mathop{\rm int}} C - C \subseteq \alpha w + B - {\mathop{\rm int}} C.$$
We conclude from  Proposition  \ref{Bxsdu} (i)  that ${\varphi _{w,B}}\left( y \right) < \alpha  = {\varphi _{w,A}}\left( y \right)$.  This completes the proof.   \hfill$\Box$

Similar to the proof of Proposition \ref{Bxxxsdu}, we can get the following proposition.
\begin{proposition}  \label{Bzxsdu}   Let $y_1, y_2 \in Y$ and $B \in {\wp _{0 - C}}\left( Y \right)$.
\begin{itemize}
\item[(i)] If  ${y_2} - {y_1} \in C$, then ${\varphi _{w,B}}\left( {{y_2}} \right) \ge {\varphi _{w,B}}\left( {{y_1}} \right)$.			
\item[(ii)] If ${y_2} - {y_1} \in {\mathop{\rm int}} C$ and $B$ is $-C$-closed, then ${\varphi _{w,B}}\left( {{y_2}} \right) > {\varphi _{w,B}}\left( {{y_1}} \right)$.			
\end{itemize}
\end{proposition}

\begin{proposition}  \label{Bxqqxsdu}   Let   $A,B, D, Q \in {\wp _{0 - C}}\left( Y \right)$.
	\begin{itemize}
		\item[(i)] If $D{ \le ^u}Q$, then ${H_w}\left( {D,B} \right) \le {H_w}\left( {Q,B} \right)$ and ${H_w}\left( {A,D} \right) \ge {H_w}\left( {A,Q} \right)$.			
		\item[(ii)] If $D{ \ll ^u}Q$,   $D$ is $-C$-compact and $B$ is $-C$-closed, then ${H_w}\left( {D,B} \right) < {H_w}\left( {Q,B} \right)$.
		\item[(iii)] If $D{ \ll ^u}Q$,   $A$ is $-C$-compact, and $D$ and $Q$ are $-C$-closed, then ${H_w}\left( {A,D} \right) > {H_w}\left( {A,Q} \right)$.		
	\end{itemize}
\end{proposition}
\emph{Proof}  (i) follows from Propositions \ref{Bxxxsdu} (i) and  \ref{Bzxsdu} (i).

(ii). By Proposition 4.7 of \cite{Han2019},  there exists  ${x_0} \in D$ such that ${H_w}\left( {D,B} \right) = {\varphi _{w,B}}\left( {{x_0}} \right)$. Due to $D{ \ll ^u}Q$, we have $D \subseteq Q - {\mathop{\rm int}} C$, and so there exists  ${y_0} \in Q$ such that ${y_0} - {x_0} \in {\mathop{\rm int}} C$. This together with  Proposition    \ref{Bzxsdu} (ii)  implies that
$${H_w}\left( {D,B} \right) = {\varphi _{w,B}}\left( {{x_0}} \right) < {\varphi _{w,B}}\left( {{y_0}} \right) \le {H_w}\left( {Q,B} \right).$$

(iii).  It follows from Proposition 4.7 of \cite{Han2019} that there exists  ${a_0} \in A$ such that ${H_w}\left( {A,Q} \right) = {\varphi _{w,Q}}\left( {{a_0}} \right)$.  Combining this with Proposition \ref{Bxxxsdu} (ii), we get
$${H_w}\left( {A,Q} \right) = {\varphi _{w,Q}}\left( {{a_0}} \right) < {\varphi _{w,D}}\left( {{a_0}} \right) \le {H_w}\left( {A,D} \right).$$
This completes the proof.   \hfill$\Box$

Similar to the proof of Lemma 4.3 of \cite{Han2022AA}, we can get the following lemma.
\begin{lemma} \label{BdcG}  If $\beta  > 0$ and $A$ is $-C$-bounded, then $A \not\subset A - \beta w - {\mathop{\rm int}} C.$
   \end{lemma}

\begin{lemma} \label{BdAAcG}  If  $A$ is nonempty  and $-C$-bounded, then ${H_w}\left( {A,A} \right) = 0$.
\end{lemma}
\emph{Proof} Noting that $A \subseteq A - C = 0w + A - C$, it follows from Proposition  \ref{Bxugsdu} that ${H_w}\left( {A,A} \right) \le 0$. Suppose that ${H_w}\left( {A,A} \right) < 0$. Then there exists $\beta  > 0$ such that ${H_w}\left( {A,A} \right) <  - \beta  $. Due to Proposition  \ref{BxYGsdu}, one has
\begin{equation}\label{DvvB1} A \subseteq \left( { - \beta } \right)w + A - {\mathop{\rm int}} C.
\end{equation}
We conclude from  Lemma \ref{BdcG}  that
$A \not\subset A - \beta w - {\mathop{\rm int}} C$, which contradicts (\ref{DvvB1}). Therefore, we have  ${H_w}\left( {A,A} \right) = 0$. This completes the proof.   \hfill$\Box$

By Theorem 5.1 of \cite{Han2019}, we can get the following lemma.
\begin{lemma} \label{lianxuuu} If $F$ is $-C$-continuous  on $K$ with nonempty  and $-C$-compact  values, then   $\gamma \left( { \cdot , \cdot } \right)$ is continuous on $K \times K$.
\end{lemma}

\begin{lemma} \label{utuuxing} Let $K$ be  a nonempty subset of $X$ and ${y _0} \in K$. Assume that   $F$ is 	strictly quasi $u$-$C$-convexlike on $K$  with nonempty, $-C$-convex and  $-C$-compact values. Then $\gamma \left( {\cdot, {y _0} } \right)$ is a 	strictly quasi  convexlike function on $K$, i.e., for any ${x_1},{x_2} \in K$ with ${x_1} \ne {x_2}$,  there exist ${x_3} \in K$ and   $\lambda \in \left[ {0,1} \right]$ such that
	$$\gamma \left( {{x_3},{y_0}} \right) < \lambda \gamma \left( {{x_1},{y_0}} \right) + \left( {1 - \lambda } \right)\gamma \left( {{x_2},{y_0}} \right).$$
\end{lemma}
{\it Proof}   Let ${x_1},{x_2} \in K$ with ${x_1} \ne {x_2}$.   Noting that $F$ is 	strictly quasi $u$-$C$-convexlike on $K$, there exist ${x_3} \in K$ and   $\lambda \in \left[ {0,1} \right]$ such that
\begin{equation}\label{kukB1} F\left( {{x_3}} \right) \subseteq \lambda F\left( {{x_1}} \right) + \left( {1 - \lambda } \right)F\left( {{x_2}} \right) - {\mathop{\rm int}} C
\end{equation}
Since  $F\left( {{x_3}} \right)$ is  $-C$-compact  and $F\left( {{y _0}} \right)$  is $-C$-closed, it follows from  Proposition 4.7 of \cite{Han2019}  that
 there exists  ${u_3} \in F\left( {{x_3}} \right)$ such that
\begin{equation}\label{kukB2}
\gamma \left( {{x_3},{y_0}} \right) = {H_w}\left( {F\left( {{x_3}} \right),F\left( {{y_0}} \right)} \right) = {\varphi _{w,F\left( {{y_0}} \right)}}\left( {{u_3}} \right).
\end{equation}
We conclude from  (\ref{kukB1})  that  there exist  ${u_1} \in F\left( {{x_1}} \right)$ and ${u_2} \in F\left( {{x_2}} \right)$  such that $\lambda {u_1} + \left( {1 - \lambda } \right){u_2} - {u_3} \in {\mathop{\rm int}} C$.  This together with  Proposition  \ref{Bzxsdu} (ii)  implies that 
\begin{equation}\label{kukB3}
{\varphi _{w,F\left( {{y_0}} \right)}}\left( {{u_3}} \right) < {\varphi _{w,F\left( {{y_0}} \right)}}\left( {\lambda {u_1} + \left( {1 - \lambda } \right){u_2}} \right).
\end{equation}
Let ${r_i} = {\varphi _{w,F\left( {{y_0}} \right)}}\left( {{u_i}} \right)$ for $i = 1,2$. For any $\varepsilon  > 0$, due to ${\varphi _{w,F\left( {{y_0}} \right)}}\left( {{u_i}} \right) < {r_i} + \varepsilon $ and Proposition  \ref{Bxsdu} (i), we have
\begin{equation}\label{kukB4}
{u_i} \in \left( {{r_i} + \varepsilon } \right)w + F\left( {{y_0}} \right) - {\mathop{\rm int}} C \subseteq \left( {{r_i} + \varepsilon } \right)w + F\left( {{y_0}} \right) - C,  \quad  i = 1,2.
\end{equation}
Since $F\left( {{y _0}} \right)$  is $-C$-convex, one has $\lambda F\left( {{y_0}} \right) + \left( {1 - \lambda } \right)F\left( {{y_0}} \right) \subseteq F\left( {{y_0}} \right) - C$. Combining this with (\ref{kukB4}), we get
\begin{eqnarray*} \lambda {u_1} + \left( {1 - \lambda } \right){u_2} &\in& \left( {\lambda {r_1} + \left( {1 - \lambda } \right){r_2} + \varepsilon } \right)w + \lambda F\left( {{y_0}} \right) + \left( {1 - \lambda } \right)F\left( {{y_0}} \right) - C \\
	& \subseteq & \left( {\lambda {r_1} + \left( {1 - \lambda } \right){r_2} + \varepsilon } \right)w + F\left( {{y_0}} \right) - C,	
\end{eqnarray*}
and so ${\varphi _{w,F\left( {{y_0}} \right)}}\left( {\lambda {u_1} + \left( {1 - \lambda } \right){u_2}} \right) \le \lambda {r_1} + \left( {1 - \lambda } \right){r_2} + \varepsilon $.
By the arbitrariness of $\varepsilon  > 0$, we have 
\begin{eqnarray*} {\varphi _{w,F\left( {{y_0}} \right)}}\left( {\lambda {u_1} + \left( {1 - \lambda } \right){u_2}} \right) &\le& \lambda {r_1} + \left( {1 - \lambda } \right){r_2} = \lambda {\varphi _{w,F\left( {{y_0}} \right)}}\left( {{u_1}} \right) + \left( {1 - \lambda } \right){\varphi _{w,F\left( {{y_0}} \right)}}\left( {{u_2}} \right) \\
&\le& \lambda {H_w}\left( {F\left( {{x_1}} \right),F\left( {{y_0}} \right)} \right) + \left( {1 - \lambda } \right) {H_w}\left( {F\left( {{x_2}} \right),F\left( {{y_0}} \right)} \right)\\
&=& \lambda \gamma \left( {{x_1},{y_0}} \right) + \left( {1 - \lambda } \right) \gamma \left( {{x_2},{y_0}} \right).	
\end{eqnarray*}
This together with   (\ref{kukB2})  and  (\ref{kukB3})  implies that 
$$\gamma \left( {{x_3},{y_0}} \right) < \lambda \gamma \left( {{x_1},{y_0}} \right) + \left( {1 - \lambda } \right)\gamma \left( {{x_2},{y_0}} \right).$$
This completes the proof.  \hfill$\Box$

Similar to the proof of Theorem  \ref{hldtx}, by Lemmas \ref{uuKhrj}, \ref{WEuuxd}, \ref{BdAAcG}, \ref{lianxuuu} and \ref{utuuxing} and Proposition \ref{BxYGsdu} and  \ref{Bxqqxsdu}  (ii), we can get the following theorem.
\begin{theorem} \label{huultx}  Let $K$ be a nonempty  and compact subset of $X$. Assume that  $K$ is    starshaped   and $F$ is  $-C$-continuous and strictly quasi $u$-$C$-convexlike  on $K$  with nonempty, $-C$-convex  and  $-C$-compact values. Then ${W_u}\left( {F,K} \right)$ and ${E_u}\left( {F,K} \right)$ are both contractible.  
\end{theorem}

\section{Contractibility of  p-minimal and  p-weak minimal solution sets}

In this section, we discuss the contractibility of  $p$-minimal and  $p$-weak minimal solution sets for set optimization problems. Han \cite{Han2020} established connectedness of weak $p$-minimal solution set for set optimization problems by using  the scalarization technique. However, he did not investigate the arcwise connectedness of $p$-minimal and  $p$-weak minimal solution sets for set optimization problems. Now, we make a new attempt to derive the arcwise connectedness of $p$-minimal and  $p$-weak minimal solution sets for set optimization problems by    utilizing strictly    cone-convexlikeness. Let ${A_C}: = \bigcap\limits_{a \in A} {\left( {a + C} \right)} $.
For $f \in B_e^*$,   let $Q\left( f \right)$
denote the set of all  $f$-solutions of (SOP), i.e.,
$$Q\left( f \right) = \left\{ {x \in K:\mathop {\sup }\limits_{a \in F\left( x \right)} f\left( a \right) \le \mathop {\inf }\limits_{\beta  \in {F}\left( y \right) _C} f\left( \beta  \right),\;\forall y \in K} \right\}.$$

By Theorem 3.1 of \cite{Han2020}, we can get the following lemma.

\begin{lemma} \label{BLHJG} If $K$ is  nonempty  and $F$ is $p$-$C$-convexlike on $K$ with nonempty, $-C$-convex  and $-C$-compact values, then $${W_p}\left( {F,K} \right) = \bigcup\limits_{f \in B_e^*} {Q\left( f \right)} .$$
\end{lemma}

\begin{remark} In Theorem 3.1 of \cite{Han2020}, we do not need to assume that  $K$ is convex. 
\end{remark}

From Step 2 in the proof of  Theorem 4.1 of \cite{Han2020}, we can get the following lemma.

\begin{lemma} \label{Qsblx} Let $K$ be a  nonempty and compact subset of $X$. If $F$ is $-C$-l.s.c on $K$, then   $Q$ is u.s.c. on $B_e^*$,  where the topology on $B_e^*$ is   the weak* topology.   
\end{lemma}

\begin{lemma} \label{Qdandj} \cite{Han2022} Let $f \in C ^*$.  If  $F$ is strictly quasi $p$-$C$-convexlike on  $K$,
	then $Q\left( f \right)$ is a   singleton.
\end{lemma}

By Remark 3.1 of \cite{Han2022} and Theorem 3.3 of \cite{Han2022}, we can get the following lemma.
\begin{lemma} \label{rkyv}   Let $K$ be a nonempty   subset of  $X$. If  $F$ is strictly quasi $p$-$C$-convexlike on $K$ with nonempty, $- C$-convex and $- C$-compact values,
	then ${W_p}\left( {F,K} \right) = {E_p}\left( {F,K} \right)$.
\end{lemma}

\begin{theorem} \label{hltx}  Let $K$ be a nonempty  and compact subset of $X$. Assume that   $F$ is $-C$-l.s.c and strictly   $p$-$C$-convexlike  on $K$  with   nonempty, $- C$-convex and $- C$-compact values. Then ${W_p}\left( {F,K} \right)$ is arcwise connected.  Moreover, ${E_p}\left( {F,K} \right)$ is  arcwise connected.
\end{theorem}
\emph{Proof}  Since  $F$ is   strictly   $p$-$C$-convexlike  on $K$, we can see that $F$ is       $p$-$C$-convexlike and strictly quasi $p$-$C$-convexlike  on $K$. 
It follows from Lemma \ref{BLHJG} that $${W_p}\left( {F,K} \right) = \bigcup\limits_{f \in B_e^*} {Q\left( f \right)} .$$ By Lemmas \ref{Qsblx} and \ref{Qdandj}, we obtain that  $Q:B_e^* \to K$ is a continuous  single-valued   mapping. It is clear that $B_e^*$ is convex, and so  $B_e^*$ is arcwise connected. Thus, we  conclude from  Lemma \ref{HltJB}  that  ${W_p}\left( {F,K} \right)$ is arcwise connected.   Combining this with  Lemma  \ref{rkyv}, we get that ${E_p}\left( {F,K} \right)$ is  arcwise connected. This completes the proof.   \hfill$\Box$

By Theorem \ref{hltx}, it is easy to get the following theorem.
\begin{theorem} \label{PXSX}  Under the assumptions of Theorem \ref{hltx},    ${W_p}\left( {F,K} \right)$ and ${E_p}\left( {F,K} \right)$ are both  contractible.
\end{theorem}

\section{An application to vector optimization problems}

In this section,  the main results presented in  Section 3  will be applied to obtain  the contractibility  of the solution sets  for vector optimization problems.
Let $f:X \to Y$ be a single valued mapping and $K$ be a nonempty subset of $X$.   We consider the following vector optimization problem (for
short, VOP): finding ${x_0} \in K$ such that
$$f\left( y \right) - f\left( {{x_0}} \right) \notin  - C\backslash \left\{ 0 \right\},\;\;\forall y \in K.$$

Let $W\left( {f,K} \right)$ denote the set of all weakly efficient solutions of (VOP), i.e.,
$$W\left( {f,K} \right) = \left\{ {x \in K:f\left( y \right) - f\left( x \right) \notin  - {\mathop{\rm int}} C,\; \forall y \in K} \right\}$$
and $E\left( {f,K} \right)$ denote the set of all efficient solutions of (VOP), i.e.,
$$E\left( {f,K} \right) = \left\{ {x \in K:f\left( y \right) - f\left( x \right) \notin  - C\backslash \left\{ 0 \right\},\;\forall y \in K} \right\}.$$

\begin{remark} \label{tuih} If $F:X \to {2^Y}$ reduces to a single-valued mapping $f:X \to Y$, then it is easy to see that ${E_l}\left( {F,K} \right) = E\left( {f,K} \right)$ and ${W_l}\left( {F,K} \right) = W\left( {f,K} \right)$.
\end{remark}

\begin{definition}  \cite{Luc}  A mapping $\phi :X \to Y$ is said to be
	$C$-continuous at ${x_0} \in X$, if
	for any neighborhood $V$ of $\phi\left( {{x_0}} \right)$, there exists
	a neighborhood $U\left( {{x_0}} \right)$ of ${x_0}$ such that
	$$\phi \left( x \right) \in V + C, \quad \forall x \in U\left( {{x_0}} \right).$$
	We say that $\phi$ is   $C$-continuous on $X$, if it is  $C$-continuous   at each point
	$x \in X$.
\end{definition}

\begin{remark} \label{JTEC} Let a  set-valued mapping $\Phi :X \to {2^Y}$ be $C$-continuous at $x_0 \in X$. If $\Phi :X \to {2^Y}$
	reduces to a single-valued mapping $\phi :X \to {Y}$, then $\phi :X \to {Y}$ is $C$-continuous and $-C$-continuous at $x_0 \in X$.	
\end{remark}

\begin{remark} \label{fzlx} It is clear that if $\phi :X \to {Y}$	is  continuous at $x_0 \in X$, then  $\phi :X \to {Y}$ is $C$-continuous and $-C$-continuous at $x_0 \in X$.	
\end{remark}

\begin{definition} \label{ygstk}  Let $D$ be a nonempty  subset of $X$. A mapping $\phi :X \to Y$ is said to be	strictly quasi $C$-convexlike on $D$, if for any ${x_1},{x_2} \in D$ with ${x_1} \ne {x_2}$,  there exist ${x_3} \in D$ and   $\lambda \in \left[ {0,1} \right]$ such that
	$$\lambda \phi \left( {{x_1}} \right) + \left( {1 - \lambda } \right)\phi \left( {{x_2}} \right) - \phi \left( {{x_3}} \right) \in {\rm{int}}C.$$
\end{definition}

From Remarks  \ref{tuih} and \ref{JTEC},   and Theorem \ref{hldtx}, we have the following theorem.
\begin{theorem} \label{ffhldtx}  Let $K$ be a nonempty  and compact subset of $X$. Assume that  $K$ is    starshaped   and $f$ is  $C$-continuous, $-C$-continuous and strictly quasi $C$-convexlike  on $K$. Then  ${W}\left( {f,K} \right)$ and  ${E}\left( {f,K} \right)$ are both contractible. 
\end{theorem}

\begin{remark}   Noting that  the  strictly quasi $C$-convexlikeness of $f$ is weaker than the strictly quasi $C$-convexity of $f$, and by Remarks \ref{Wddfk} and \ref{fzlx},   we can see that Theorem \ref{ffhldtx}  improves  Corollary 4.14 of \cite{Luc}.
\end{remark}

\vskip 6mm
\noindent{\bf Acknowledgments}

\noindent  This work was supported by the National Natural Science Foundation of China [grant number 11801257], [grant number 12071165], the Natural Science Foundation of Jiangxi Province [grant number 20232BAB211012].

\end{document}